\documentclass{amsproc}

\usepackage{a4wide}
\usepackage{amsmath}
\usepackage{enumerate}

\usepackage{amsmath,amsthm}
\usepackage{amsfonts}
\usepackage{amssymb}
\usepackage{graphicx, color, soul}
\usepackage[T1]{fontenc}
\usepackage{verbatim}
\newtheorem{rem}{Remark}{\bf}{\rm}{\rm}

\newtheorem*{theorem*}{Theorem}
\newtheorem*{cor*}{Corollary}

\newtheorem{theorem}{Theorem}
\newtheorem{cor}{Corollary}
\newtheorem{prop}{Proposition}

\newtheorem{lem}{Lemma}

{\rm}{\rm}

\newcommand{\Hol}{\operatorname{Hol}}
\newcommand{\GL}{\operatorname{GL}}
\newcommand{\SO}{\operatorname{SO}}
\newcommand{\Un}{\operatorname{U}}

\newcommand{\Spin}{\operatorname{Spin}}

\newcommand{\pr}{\operatorname{pr}}
\newcommand{\id}{\operatorname{id}}
\newcommand{\rk}{\operatorname{rk}}

\def\so{\mathfrak{so}}
\def\su{\mathfrak{su}}
\def\u{\mathfrak{u}}
\def\t{\mathfrak{t}}
\def\ff{\mathfrak{f}}

\newcommand{\h}{\mathfrak{h}}
\newcommand{\hol}{\mathfrak{hol}}
\newcommand{\gl}{\mathfrak{gl}}

\newcommand{\Real}{\mathbb{R}}
\newcommand{\Co}{\mathbb{C}}

\begin{document}
	
	\title{Holonomy of K-contact sub-Riemannian manifolds}
	

	\author{Anton S. Galaev}  
		\thanks{$^1$
		University of Hradec Kr\'alov\'e, Faculty of Science, Department of Mathematics, Rokitansk\'eho 62, 500~03 Hradec Kr\'alov\'e,  Czech Republic
		}

	\maketitle

	\begin{center}
		In memory of my father Sergey V. Galaev 
	\end{center}

	\begin{abstract}
		
It is shown that the horizontal holonomy group of a K-contact sub-Riemannian manifold either coincides with the holonomy group of a Riemannian manifold, or it is a codimension-one normal subgroup  of the later group. This implies the classification of the holonomy algebras of K-contact sub-Riemannian manifolds.  
The question of existence of parallel horizontal spinors, examples, and consequences are discussed. 



		\vskip0.1cm
		
		\vskip0.1cm
		
		
		
	\end{abstract}

\section{Introduction}

Let $D$ be a distribution on a smooth manifold $M$. A horizontal  connection $\nabla$ defines parallel transport of vectors tangent to $D$ along curves tangent to $D$. This gives rise to the so-called horizontal holonomy. For the first time the horizontal holonomy was investigated in the contact sub-Riemannian case \cite{FGR}. Later this notion was generalized and studied in several works, e.g., in \cite{CGJK,Gro,HMCK,Mal}.  

In the present paper we study holonomy of K-contact sub-Riemannian manifolds $(M,\theta,g)$. Here $M$ is a smooth manifold with a contact form $\theta$, and $g$ is a sub-Riemannian metric on the contact distribution $D=\ker\theta$. We assume that the Reeb vector field $\xi$ is complete and its flow   preserves the metric $g$. 

In Section \ref{secHorCon}, we recall the definition of Schouten connection $\nabla$, which is the horizontal Levi-Civita connection on a sub-Riemannian manifold $(M,D,g)$. In Section \ref{secExtCon}, we discuss the extensions of the horizontal connection $\nabla$ to connections on the vector bundle $D$. Among these extensions are the adapted connection $\nabla^0$ and the Wagner connection $\nabla^W$. The connections $\nabla^W$ and $\nabla^0$ may be extended to metric connections with torsion on $M$, and their holonomy may be interesting also from this point of view, see \cite{A10}. 

In Section \ref{secHorHol}, we recall some facts about the horizontal holonomy from \cite{FGR,HMCK}. In particular, the holonomy of the connections $\nabla$ and $\nabla^W$ coincide, while the holonomy of the connection $\nabla$ is contained in the holonomy of the connection $\nabla^0$.
The Reeb vector field $\xi$ defines a one-dimensional foliation $\mathcal{F}$ on $M$ (the characteristic foliation), by this reason, in Section \ref{secTransStrFol}, we discuss transverse structures on foliations. In Section \ref{secConTheta}, we compute the holonomy of the connection in the line bundle 
$\left<\xi\right>$ defined by the contact form $\theta$.

In Section \ref{secversus}, we compare the holonomy groups of the connections $\nabla$ and $\nabla^0$. We prove the main result of the paper stating that either these groups coincide, or the holonomy group of  $\nabla$ is a codimension-one subgroup of the holonomy group of $\nabla^0$. 
In Section \ref{secClassif}, we study the second possibility in the case of simply connected spaces. We show that this possibility is equivalent to the existence of a transverse K\"ahler structure on the characteristic foliation and a relation between the 2-form $d\theta$ and the transverse Ricci form.

Section \ref{secEx} deals with examples and consequences of the results of Section \ref{secClassif}. We give a construction illustrating the results from Section \ref{secClassif}, we consider the special case of Sasaki manifolds, and we discuss the relation of the holonomy of the connection $\nabla$ to the existence of parallel horizontal spinors.

\vskip0.3cm

{\bf Acknowledgements.}
The author is thankful to D.V. Alekseevsky, Ya. Bazaykin,  I.~Chrysikos for useful and stimulating discussions, and to anonymous referee for careful reading of the paper and pointing out several shortcomings.
The research was supported by the project GF24-10031K of Czech Science Foundation (GA\v{C}R).

\section{Horizontal connection}\label{secHorCon}

Let $D$ be a distribution on a smooth manifold $M$. A piece-wise smooth curve in $M$ is called {\it horizontal} if it is tangent to $D$.
{\it A horizontal  connection} on the distribution $D$ is a map 
\begin{equation}\label{nablaD}\nabla:\Gamma(D)\times \Gamma(D)\to \Gamma(D)\end{equation} sharing the usual properties of the covariant derivative. For a given horizontal connection $\nabla$, a horizontal curve $\gamma:[a,b]\to M$ defines the parallel transport
\begin{equation}\label{tau}\tau_\gamma: D_{\gamma(a)}\to D_{\gamma(b)}.\end{equation} along $\gamma$.
Let $g$ be a sub-Riemannian metric, i.e., a  field of non-degenerate positively definite symmetric bilinear forms on the distribution $D$. Suppose that a rigging $D'$ is fixed, i.e., there is a decomposition
$$TM=D\oplus D'$$ into the direct sum of the distributions. 
Let $$\pi:\Gamma(TM)\to\Gamma(D),\quad \pi':\Gamma(TM)\to\Gamma(D')$$ be the corresponding projections.
The metric $g$ and the rigging $D'$ define the Schouten connection that is the horizontal Levi-Civita connection given by the Koszul formula
\begin{equation}\label{SchoutenCon}
2g(\nabla_X Y,Z)= Xg(Y,Z) +Yg(Z,X) - Zg(X,Y)+ g(\pi[X,Y],Z)- g(\pi[Y,Z],X)- g(\pi[X,Z],Y),
\end{equation} where $X,Y,Z\in\Gamma(D)$.
This is the unique connection preserving the metric $g$ and having zero torsion in the sense that 
\begin{equation}\label{torsioncondit}\nabla_X Y-\nabla_Y X-\pi[X,Y]=0,\quad\forall X,Y\in\Gamma(D).\end{equation}
The curvature tensor of a connection $\nabla$ was defined by Schouten
in the following way:
\begin{equation}\label{Schouten}
R(X,Y)Z=\nabla_X\nabla_Y Z-\nabla_Y\nabla_X Z-\nabla_{\pi[X,Y]}Z-\pi[\pi'[X,Y],Z],\quad\forall X,Y,Z\in\Gamma(D).
\end{equation}
The Schouten tensor satisfies the Bianchi identity 
\begin{equation}\label{BianchiSch}
R(X,Y)Z+R(Y,Z)X+R(Z,X)Y=0,  \quad\forall X,Y,Z\in\Gamma(D).
\end{equation}

\section{Extended connections}\label{secExtCon}

Let, as in the previous section, $\nabla$ be a horizontal connection on a distribution $D$. {\it An extended connection} is a connection $$\nabla^D:\Gamma(TM)\times \Gamma(D)\to \Gamma(D)$$ on the vector bundle $D$ such that $$\nabla^D_XY=\nabla_XY,\quad  \forall X,Y\in\Gamma(D).$$   
Wagner \cite{Wagner37} noted that the flatness of a horizontal connection is generally not equivalent to the vanishing of its Schouten tensor, and he constructed a special extended connection such that its curvature vanishes if and only if the horizontal connection is flat. The extended connection and its curvature tensor are called nowadays the Wagner connection and the Wagner curvature tensor,  see \cite{DG} for the modern exposition. 

Now we start to consider the case of contact sub-Riemannian manifolds.
Assume that the distribution $D$ coincides with the kernel of a fixed contact form $\theta\in\Omega^1(M)$. Denote by $\xi$ the corresponding Reeb vector field on $M$ defined by the conditions 
$$\theta(\xi)=1,\quad \iota_X d\theta=0.$$  
We assume that $\dim M=2m+1\geq 5$. 

We denote by the triple $(M,\theta,g)$ 
a contact manifold $(M,\theta)$ with a sub-Riemannian metric $g$ on the contact distribution and call it {\it a contact sub-Riemannian manifold}.

An arbitrary extended connection is uniquely defined by a field $N$ of endomorphisms of $D$ and satisfies the conditions
\begin{equation}\label{nablaN}
\nabla^N_\xi X=[\xi,X]+NX,\quad \nabla^N_XY=\nabla_XY,\quad\forall X,Y\in\Gamma(D).
\end{equation} 

The curvature tensor of the connection $\nabla^N$ is given by
\begin{align} R^N(X,Y)&=R(X,Y)+d\theta(X,Y)N,\quad\forall X,Y\in\Gamma(D)\\
R^N(\xi,X)&=(\mathcal{L}_\xi\nabla)_X-\nabla_XN,\quad\forall X\in\Gamma(D),
\end{align}
see, e.g., \cite{Papa2016,Papa2015}. Here $R$ is the Schouten tensor defined in the previous section.

The extended connection given by $N$ satisfying
$$g(NX,Y)=\frac{1}{2}(\mathcal{L}_\xi g)(X,Y),\quad \forall X,Y\in\Gamma(D),$$ is called {\it  the adapted connection}, see, e.g., \cite{FGR}. 

For vector fields $X,Y\in\Gamma(D)$ we define the bivector
$X\wedge Y$ as the tensor field 
$$X\wedge Y=X\otimes Y-Y\otimes X.$$ The module of all bivectors is denoted by $\Gamma(\wedge^2 D)$. For a 2-form 
$\omega$ on $D$ with values in a vector bundle over $M$, we set $$\omega(X\wedge Y)=2\omega(X,Y).$$
This allows us to consider the values $\omega(\alpha)$ of the 2-form $\omega$ on any bivector $\alpha\in \Gamma(\wedge^2 D)$.

The Wagner connection $\nabla^W$ in the case of a contact sub-Riemannian manifold is defined in the following way, see, e.g., \cite{Papa2016,Papa2015}.
Denote by $(d\theta)^{-1}\in\Gamma(\wedge^2 D)$ the bivector  such that its matrix of coordinates is the inverse to the matrix of coordinates of the form $d\theta$ with respect to any local frame of $D$. It holds \begin{equation}\label{dthetadtheta-1}d\theta((d\theta)^{-1})=-4m.\end{equation}
The Wagner connection \cite{Wagner37,Wagner41} is defined by the condition on its curvature tensor
\begin{equation}\label{Wagnercond}
R^W((d\theta)^{-1})=0.\end{equation}
The corresponding field $N$ is given by
\begin{equation}\label{NWagner} N=\frac{1}{4m}R((d\theta)^{-1}).\end{equation}
 This extended connection is a particular  {\it  special extension} defined in \cite{FGR}. For the case of CR-manifolds such connection was constructed in \cite{Leitner18}, where it was called {a basic connection}.

A contact sub-Riemannian manifold is called {\it K-contact} if $\mathcal{L}_\xi g=0$. It is clear that in this case the corresponding   adapted connection is given by the choice $N=0$. We denote this connection by $\nabla^0$. It holds 
\begin{equation}\label{nablaN0}
\nabla^0_\xi X=[\xi,X],\quad \nabla^0_XY=\nabla_XY,\quad\forall X,Y\in\Gamma(D).
\end{equation} 
 The curvature tensor of this connection is given by 
\begin{equation}\label{R0} R^0(X,Y)=R(X,Y),\quad R^0(\xi,X)=0,\quad \forall X,Y\in \Gamma(D).\end{equation}
 The curvature tensor of the Wagner connection in the K-contact case satisfies  
\begin{equation}\label{RW} R^W(X,Y)=R(X,Y)+d\theta(X,Y)N,\quad R^W(\xi,X)=-\nabla_XN,\quad \forall X,Y\in \Gamma(D),\end{equation}
where $N$ is given by \eqref{NWagner}.


\section{Horizontal holonomy}\label{secHorHol}

Given a distribution $D$ on a manifold $M$, a horizontal connection $\nabla$ as in \eqref{nablaD}, and a point $x\in M$, one defines the {\it (horizontal) holonomy group} of the connection $\nabla$ as the group of parallel transports along piece-wise smooth horizontal loops at the point $x\in M$. We denote this group by $\Hol_x(\nabla)$. The {\it restricted holonomy group} $\Hol^0_x(\nabla)$ is the subgroup of 
$\Hol_x(\nabla)$ corresponding to contactable loops. Note that if the distribution $D$ is bracket-generating, then any contractable horizontal loop is contactable in the class of horizontal loops, see, e.g., \cite{CGJK}. The holonomy group $\Hol_x(\nabla)$ is a Lie subgroup of $\GL(D_x)$, and the restricted holonomy group $\Hol^0_x(\nabla)$ is the identity component of $\Hol_x(\nabla)$, \cite{CGJK,FGR,HMCK}.
The corresponding Lie subalgebra of $\gl(D_x)$ is called {\it the holonomy algebra} and is denoted by $\hol_x(\nabla)$. The parallel transport along a horizontal curve defines the isomorphism of the holonomy at the initial and the end-point of the curve, so the holonomy is independent on the point.

Consider now a contact sub-Riemannian manifold $(M,\theta,g)$ as in Section \ref{secExtCon}.
The following are versions of the Ambrose-Singer holonomy Theorem for the horizontal connection \cite{FGR}, see also \cite{CGJK}.

\begin{theorem} The Lie algebra $\hol_x(\nabla)$ is spanned by 
	the following endomorphisms of $D_x$:
	$$\tau_\gamma^{-1}\circ R(\alpha)\circ\tau_\gamma,$$
	where $\gamma$ is an arbitrary piece-wise smooth horizontal curve starting at $x$, and $\alpha$ is a bivector from $D$ at the end-point of $\gamma$ such that $d\theta(\alpha)=0$. 
\end{theorem}

\begin{theorem}\label{thASnabla} The Lie algebra $\hol_x(\nabla)$ is spanned by 
	the following endomorphisms of $D_x$:
	$$\tau_\gamma^{-1}\circ R^W(X,Y)\circ\tau_\gamma,$$
	where $\gamma$ is an arbitrary piece-wise smooth horizontal curve starting at $x$, and $X$, $Y$  are vectors from $D$ at the end-point of $\gamma$. 
\end{theorem}

In \cite{CGJK}, was introduced a notion of {\it selector} for a class of distributions. For a contact distribution $D$ defined by the kernel of a 1-form $\theta$, a selector is a bivector $\alpha\in\Gamma(\wedge^2 D)$ such that the function $d\theta(\alpha)$ is non-vanishing. From \cite[Theorem 2.16]{CGJK} it follows

\begin{theorem} Let $(M,\theta,g)$ be a contact sub-Riemannian manifold and let $\nabla^N$ be an extension of the Schouten connection $\nabla$. Suppose that there is a selector $\alpha$ of $D$ such that $R^N(\alpha)=0$, then $\Hol(\nabla)=\Hol(\nabla^N)$.  \end{theorem}

This theorem, \eqref{dthetadtheta-1}, and \eqref{Wagnercond} imply

\begin{theorem}\label{thholW} Let $(M,\theta,g)$ be a contact sub-Riemannian manifold. Then the horizontal holonomy group $\Hol(\nabla)$ coincides with the holonomy group $\Hol(\nabla^W)$ of the Wagner connection. 
\end{theorem}

\section{Transverse structures of  foliations}\label{secTransStrFol}

In this section we follow the exposition from \cite{HV}.
Let $\mathcal{F}$ be a foliation on a smooth manifold $M$. Denote by  $L\subset TM$ the distribution  tangent  to $\mathcal{F}$. Let $Q=TM/L$ be the normal bundle of $\mathcal{F}$. {\it A transverse tensor field} on   $\mathcal{F}$ is a tensor field $A$  on the bundle $Q$ satisfying the condition
$$\mathcal{L}_\xi A=0,\quad\forall \xi \in\Gamma(L).$$
Using an auxiliary Riemannian metric on $M$, one may identify $Q$ with a distribution on $M$ such that it holds
$$TM=L\oplus Q.$$ Consider the corresponding projections
$\pi:TM\to Q$ and $\pi':TM\to L$. Equality \eqref{SchoutenCon} with $X,Y,Z\in\Gamma(Q)$ defines the horizontal Levi-Civita connection $\nabla$ on $Q$. {\it The transversal Levi-Civita connection} is the connection on the vector bundle $Q$ defined by
$$\nabla^T_XY=\nabla_XY,\quad \nabla^T_\xi Y=\pi[\xi,Y],\quad X,Y\in\Gamma(Q),\quad\xi\in\Gamma(L).$$
The holonomy of the connection $\nabla^T$ is called {\it the transverse} holonomy. It is clear that a tensor field $A$ defined on $Q$ is transverse if and only if
$$\nabla^T_\xi A=0,\quad \forall\xi\in\Gamma(L).$$ In particular, each $\nabla^T$-parallel tensor field on $Q$ is a transverse tensor field.

{\it A transverse K\"ahler structure} on  $\mathcal{F}$ is a transverse Hermitian structure $J$ on $(Q,g)$ parallel  
with respect to the transversal connection $\nabla^T$.
A transverse K\"ahler structure exists if and only if $\Hol(\nabla^T)$ is contained in $\Un(m)$, where $\rk Q=2m$. The Ricci form of a transverse K\"ahler structure is defined as $$\rho^T(X,Y)=\mathrm{Ric}^T(JX,Y),\quad X,Y\in\Gamma(Q).$$
Here $\mathrm{Ric}^T$ is the transverse Ricci tensor 
defined by $$\mathrm{Ric}^T(X,Y)=\sum_{k=1}^{2m} g (R^T(e_k,X)Y,e_k),$$
where $e_1,\dots, e_{2m}$ is an arbitrary local $g$-orthonormal frame on $Q$.
The Ricci form $\rho^T$ coincides with the induced connection 
on the bundle $K$ of complex $m$-forms. For a transversal K\"ahler structure, the equality $\rho^T=0$ is equivalent to the condition $\hol(\nabla^T)\subset\su(m)$.

Let $(M,\theta,g)$ be a K-contact sub-Riemannian manifold. Then the Reeb vector field $\xi$ defines a one-dimensional foliation $\mathcal{F}$. The contact distribution $D$ may be identified  with the normal bundle of $\mathcal{F}$, the metric $g$ is a transverse metric on $\mathcal{F}$, and the adapted connection $\nabla^0$ coincides with the transversal connection~$\nabla^T$.  

\section{The connection $\nabla^\theta$}\label{secConTheta}

Let $(M,\theta)$ be a contact manifold.  
Using the 1-form $\theta$, we define the connection $\nabla^\theta$ on the trivial line bundle $\left<\xi\right>$ over $M$ by setting 
$$\nabla^\theta_X\xi=\theta(X)\xi,\quad X\in\Gamma(TM).$$
It is clear that the curvature of the connection $\nabla^\theta$ coincides with $d\theta$, i.e.,
$$R^\theta(X,Y)\xi=(d\theta(X,Y))\xi,\quad X,Y\in\Gamma(TM).$$
Let us compute the parallel transport for the connection $\nabla^\theta$. Let $\mu(t)$, $t\in [a,b]$, be a piece-wise smooth curve in $M$. Let $\lambda(t)\xi_{\gamma(t)}$ be a section of the bundle $\left<\xi\right>$ along the curve $\gamma(t)$. The parallel transport equation takes the form
$$\dot\lambda(t)+\lambda(t)\theta(\dot\mu(t))=0.$$
The solution of this equation is
$$\lambda(t)=\lambda(a)\exp\left(-\int_a^t\theta(\dot\mu(s))ds\right).$$
Consequently, if the curve $\mu(t)$ is closed, then the parallel transport along $\mu(t)$ is given by
$$\lambda(b)=\lambda(a)\exp\left(-\int_\mu \theta\right).$$
This implies that the holonomy group $\Hol_x(\nabla^\theta)$ is isomorphic to the subgroup $\Real_+\subset\GL(1,\Real)$.

\section{$\Hol(\nabla)$ versus $\Hol(\nabla^0)$}\label{secversus}

Let $(M,\theta,g)$ be a K-contact sub-Riemannian manifold. Denote by $\nabla^0$ the extended connection given, as in Section \ref{secExtCon}, by $N=0$. We assume that the Reeb vector field $\xi$ is complete and we denote its flow by $\varphi_t$. In this section we compare the holonomy groups $\Hol_x(\nabla)$ and $\Hol_x(\nabla^0)$  at  a fixed point $x\in M$. 

If the orbit $\varphi_t(x)$ of the point $x$ is cyclic, then we denote by $t_0$ the smallest positive $t_0$ such that $\varphi_{t_0}(x)=x$. Otherwise we assume that $t_0=0$.  

Let $$\mu:[a,b]\to M$$ be a piece-wise curve. We define the new curve   
\begin{equation}\label{deftildmu}
\tilde \mu(t)=\varphi_{f(t)}\mu(t),\quad t\in[a,b],\end{equation}
where 
$$f(t)=-\int_a^t\theta(\dot\mu(s))ds.$$

\begin{lem} The curve $\tilde \mu$ is horizontal. 
\end{lem}

{\bf Proof}. Differentiating \eqref{deftildmu}, we get
$$\dot{\tilde\mu}(t)=\dot f(t) \xi_{\varphi_{f(t)}\mu(t)}+(d\varphi_{f(t)})\dot\mu(t)=-\theta(\dot\mu(t))\xi_{\varphi_{f(t)}\mu(t)}+(d\varphi_{f(t)})\dot\mu(t),$$ where 
$$d\varphi_{f(t)}:T_{\mu(t)}M\to T_{\varphi_{f(t)}\mu(t)}M$$ is the differential of the diffeomorphism $\varphi_{f(t)}$. Since the flow $\varphi$ preserves the vector field $\xi$, it holds 
$$\xi_{\varphi_{f(t)}\mu(t)}=d\varphi_{f(t)}\xi_{\mu(t)}.$$
This implies that $$\dot{\tilde\mu}(t)=(d\varphi_{f(t)})\left(\dot\mu(t)-\theta(\dot\mu(t))\xi_{\mu(t)}\right).$$
Since $\theta(\xi)=1$, the vector $\dot\mu(t)-\theta(\dot\mu(t))\xi_{\mu(t)}$ is horizontal. Since $\varphi_t$ preserves the contact distribution~$D$, we get that the vector $\dot{\tilde\mu}(t)$ is horizontal. \qed  

Since the flow $\varphi_t$ consists of automorphism of the distribution $D$, and $L_\xi g=0$, we see that the flow $\varphi_t$ consists of automorphisms of the connection $\nabla$. We will need also the following statement.

\begin{lem}\label{lemautom}
	The flow $\varphi_t$ consists of automorphism of the connection $\nabla^0$.
\end{lem}

{\bf Proof.} It is enough to prove that $(d\varphi_t)\nabla^0_\xi X=\nabla^0_{(d\varphi_t)\xi} (d\varphi_t)X$ for all $X\in\Gamma(D)$. We have
$$(d\varphi_t)\nabla^0_\xi X=(d\varphi_t)[\xi, X]=[(d\varphi_t)\xi, (d\varphi_t)X]=[\xi, (d\varphi_t)X]=\nabla^0_{\xi} (d\varphi_t)X=
\nabla^0_{(d\varphi_t)\xi} (d\varphi_t)X.$$
\qed

\begin{lem} Let $X(t)$ be a $\nabla^0$-parallel horizontal vector field along a curve $\mu(t)$.
	Then the  vector field
	$$\tilde X(t)=(d\varphi_{f(t)})X(t)$$
	is horizontal and it is $\nabla$-parallel along the curve $\tilde\mu(t)$. 
\end{lem}

{\bf Proof.} Consider the subset
$$W\subset\Real^2$$
consisting of the pairs $(t,r)\in\Real^2$ such that $t\in [a,b]$, and
$$r\in\left\{ \begin{matrix}[0,f(t)], \text{if } f(t)\geq 0,\\ [f(t),0], \text{if } f(t)<0.   \end{matrix}  \right.$$
Consider the map $$F:W\to M$$ defined by
$$F(t,r)=\varphi_r\mu(t).$$
Define the vector field $$Y(t,r)= (d\varphi_{r})X(t)$$ along $F$.
It is obvious that $$\mu(t)=F(t,0), \quad \tilde\mu(t)=F(t,f(t)),$$ $$X(t)=Y(t,0), \quad \tilde X(t)=Y(t,f(t)).$$ 

From Lemma \ref{lemautom} and the definition of the vector field $Y$ it follows that $$\nabla^0_{\xi_{F(t,r)}}Y(t,r)=0.$$ In particular,
$$\nabla^0_{\xi_{\mu(t)}}Y(t,0)=0.$$


Finally, \begin{multline*}\nabla_{\dot{\tilde\mu}(t) } \tilde X(t)= \nabla_{(d\varphi_{f(t)})\left(\dot\mu(t)-\theta(\dot\mu(t))\xi_{\mu(t)}\right)}(d\varphi_{f(t)})X(t)=\nabla_{\dot\mu(t)-\theta(\dot\mu(t))\xi_{\mu(t)}} X(t)\\=\nabla^0_{\dot\mu(t)-\theta(\dot\mu(t))\xi_{\mu(t)}} Y(t,0)=\nabla^0_{\dot\mu(t)}Y(t,0)-\theta(\dot\mu(t))\nabla^0_{\xi_{\mu(t)}}Y(t,0)= \nabla^0_{\dot\mu(t)}X(t) =0.\end{multline*} \qed

\begin{cor}\label{corthesametrans} Let $\mu:[a,b]\to M$ be a piece-wise smooth curve in $M$, and let $\tilde\mu:[a,b]\to M$ be the corresponding horizontal curve. Suppose that $\tilde\mu(b)=\mu(b)$. Then the curves $\tilde \mu$ and $\mu$ are homotopic, and  the parallel transport of the connection $\nabla^0$ along the curve $\mu$ coincides with the parallel transport of the connection $\nabla$ along the curve $\tilde\mu$, i.e., $$\tau^0_\mu=\tau_{\tilde\mu}.$$ \end{cor}

\begin{lem}\label{lemalsoloop} Suppose that $\mu$ is a loop at the point $x$. Then the curve $\tilde\mu$ is also a loop at the point $x$  if and only if
	\begin{equation}\label{condintegr}\int_\mu\theta\in \mathbb{Z}t_0.\end{equation}
\end{lem}

{\bf Proof.} Let $\mu:[a,b]\to M$ be a loop at the point $x$. By the definition of the curve $\tilde\mu$, it holds
$$\tilde\mu(b)=\varphi_cx,$$ where
$$c=-\int_\mu\theta.$$ From this and the definition of $t_0$ it follows that the condition $\tilde\mu(b)=x$ is equivalent to the condition $c\in\mathbb{Z}t_0$. \qed
 
 From 
 Corollary \ref{corthesametrans} and Lemma \ref{lemalsoloop} it follows

\begin{cor}\label{coralsoloop} If $\mu$ is a loop at the point $x$ satisfying \eqref{condintegr}, then $\tau^0_\mu\in \Hol_x(\nabla)$.
\end{cor}

Corollary \ref{coralsoloop} and Lemma \ref{lemalsoloop} imply

\begin{prop} An element $h\in\Hol_x(\nabla^0)$ may be represented by a loop satisfying \eqref{condintegr} if and only if $h\in\Hol_x(\nabla)$, i.e.,
	it holds $$\Hol_x(\nabla)=\left\{h\in\Hol_x(\nabla^0)\left|\exists \mu: \tau^0_\mu=h, \int_\mu\theta\in\mathbb{Z}t_0\right\}\right..$$
\end{prop}

\begin{prop}
	The Lie subgroup  $\Hol_x(\nabla)\subset \Hol_x(\nabla^{0})$ is normal.
\end{prop}

{\bf Proof.} Let $\mu$ be a horizontal loop at the point $x$, and let $\gamma$ be an arbitrary loop at $x$. Then $$\int_{\gamma\mu\gamma^{-1}}\theta=\int_{\gamma}\theta+\int_\mu\theta+\int_{\gamma^{-1}}\theta=\int_\mu\theta=0,$$ which implies that
$$(\tau^{0}_{\gamma}	)^{-1}\tau_\mu\tau^{0}_\gamma=(\tau^{0}_{\gamma}	)^{-1}\tau^{0}_\mu\tau^{0}_\gamma=\tau^{0}_{\gamma\mu\gamma^{-1}}\in \Hol_x(\nabla).$$ \qed

Let $\gamma$ and $\mu$ be two loops at the point $x$ satisfying the condition
$$ \int_\mu\theta=\int_\gamma\theta+kt_0,\quad k\in\mathbb{Z}.$$
Then it holds that $\tau^{0}_{\mu\gamma^{-1}}\in\Hol_x(\nabla)$. This shows that the maps
$$\lambda:\Hol_x(\nabla^\theta)\to\Hol_x(\nabla^{0})/\Hol_x(\nabla),\quad \lambda':\Hol_x(\nabla^\theta)/e^{\mathbb{Z}t_0}\to\Hol_x(\nabla^{0})/\Hol_x(\nabla),$$
$$\lambda:\exp\left(-\int_\mu\theta\right)\mapsto \tau^{0}_\mu\cdot \Hol_x(\nabla),\quad \lambda':
\exp\left(-\int_\mu\theta\right)\,\,\cdot e^{\mathbb{Z}t_0} \mapsto \tau^{0}_\mu\cdot \Hol_x(\nabla)$$
are well-defined. It is obvious that the both maps are surjective  Lie group homomorphisms.

We have proved the following theorem.

\begin{theorem}\label{thHolnabl}
	Let $(M,\theta,g)$ be a K-contact sub-Riemannian manifold. Suppose that the Reeb vector field $\xi$ is complete. Then one of following conditions holds:
	\begin{itemize}
		\item[1.] $\Hol_x(\nabla)=\Hol_x(\nabla^{0})$;
		\item[2.]  $\Hol_x(\nabla)\subset\Hol_x(\nabla^{0})$ is a normal subgroup of codimension one.
				\end{itemize} 
If the  condition 2 holds true, then  the factor group $\Hol_x(\nabla^{0})/\Hol_x(\nabla)$ is connected, and it holds 
\begin{equation}\label{Hol/Hol0}\Hol_x(\nabla)/\Hol^0_x(\nabla)= \Hol_x(\nabla^{0})/\Hol^0_x(\nabla^{0}).\end{equation}
\end{theorem}


Suppose that $\Hol_x(\nabla)\subset\Hol_x(\nabla^{0})$ is a normal Lie subgroup of codimension one, i.e., 
$\hol_x(\nabla)\subset\hol_x(\nabla^{0})$ is  an ideal of codimension one. Since the subalgebra $\hol_x(\nabla^{0})\subset \so(2m)$ is compact, there exists a complementary one-dimensional ideal $\t_x$, i.e., it holds 
$$\hol_x(\nabla^{0})=\hol_x(\nabla)\oplus\t_x.$$
The holonomy algebra $\hol_x(\nabla^{0})\subset\so(2m)$ may be decomposed as the direct sum
$$\hol_x(\nabla^{0})=\hol_x(\nabla^{0})'\oplus \mathfrak{z}(\hol_x(\nabla^{0}))$$
of the commutator $\hol_x(\nabla^{0})'$ (the semisimple part of $\hol_x(\nabla^{0})$) and the center $\mathfrak{z}(\hol_x(\nabla^{0}))$. It is obvious that
$$\t_x\subset \mathfrak{z}(\hol_x(\nabla^{0})).$$ Let $$\t_x^\bot\subset \mathfrak{z}(\hol_x(\nabla^{0}))$$ be the orthogonal complement to $\t_x$ in $\mathfrak{z}(\hol_x(\nabla^{0}))$ 
with respect to  the extension of $g_x$ to $\so(2m)$ (which is proportional to the Cartan-Killing form of $\so(2m)$). We thus conclude that
$$\hol_x(\nabla)=\hol_x(\nabla^{0})'\oplus\t_x^\bot.$$

 Let $T_x\subset\Hol^0_x(\nabla^{0})$ be the connected normal Lie subgroup corresponding to the ideal $\t_x\subset\hol_x(\nabla^{0})$, then it holds
$$\Hol^0_x(\nabla^{0})=\Hol^0_x(\nabla)\cdot T_x.$$ This and \eqref{Hol/Hol0} imply that 
$$\Hol_x(\nabla^{0})=\Hol_x(\nabla)\cdot T_x.$$

\section{Classification}\label{secClassif} 

Let, as above, $(M,\theta,g)$ be a K-contact sub-Riemannian manifold with the complete Reeb vector field $\xi$. 

\begin{prop}\label{propholnabl0} The connected Lie group $\Hol^0_x(\nabla^0)\subset\SO(D_x,g_x)=\SO(2m)$ is the holonomy group of the Levi-Civita connection of a Riemannian manifold.
\end{prop}

{\bf Proof.} The connection $\nabla^0$ may be extended to a linear connection on $M$ by setting $\nabla^0\xi=0$. It is clear that the curvature and holonomy of this extension coincide with the initial $R^0$ and $\Hol_x(\nabla^0)$, respectively.  
From the classical Ambrose-Singer holonomy Theorem for linear connections, \eqref{BianchiSch}, and \eqref{R0} it follows that $\hol_x(\nabla^{0})\subset\so(2m)$ is a Berger subalgebra, which means that it is spanned by the images of algebraic curvature tensors \cite{Besse}. Each Berger subalgebra of $\so(2m)$ is the holonomy algebra of a Riemannian manifold \cite{Besse}. This implies the statement of the proposition. \qed 

In the rest of this section we will be interested in the holonomy algebra $\hol_x(\nabla)$, so we will assume that the manifold $M$ is simply connected.

By Proposition \ref{propholnabl0}, $\hol_x(\nabla^0)\subset\so(D_x)$ is the  holonomy algebra of a Riemannian manifold. 
According to the de~Rham decomposition Theorem,  there exists an $\hol_x(\nabla^{0})$-invariant orthogonal decomposition of the vector space $D_x$,
\begin{equation}\label{deRD} D_x=D^0_x\oplus D_x^1\oplus\cdots\oplus D_x^r\end{equation}
and the decomposition of $\hol_x(\nabla^{0})$ into the direct sum of ideals \begin{equation}\label{deRh}\hol_x(\nabla^{0})=\h_1\oplus\cdots\oplus \h_r\end{equation}
such that $\h_i\subset \so(D_x^i)$ is an irreducible Riemannian holonomy algebra for $i=1,\dots, r$, and  $\hol_x(\nabla^{0})$ annihilates $D^0_x$.
According to the holonomy principle, the $\hol_x(\nabla^0)$-invariant subspaces $D_x^i\subset D_x$ define  $\nabla^{0}$-parallel distributions $D^i$ over $M$. 

Fix an $i\in\{0,\dots,r\}$ and consider the distribution
$$L_i=\left<\xi\right>\oplus D^0\oplus D^1\oplus\cdots\oplus D^{i-1}\oplus D^{i+1}\oplus\cdots\oplus D^r.$$ 
It holds that $$TM=L_i\oplus D^i.$$
From the fact that the distributions $D^j$ are $\nabla^0$-parallel, \eqref{torsioncondit}, and \eqref{nablaN0} it follows that the distribution $L_i$ is involutive. Denote by $\mathcal{F}_i$ the corresponding foliation on $M$. The distribution $D^i$ may be identified with the transversal distribution to $\mathcal{F}_i$. From the fact that the distributions $D^j$ are $\nabla^0$-parallel and \eqref{torsioncondit} it follows that the restriction $g_{D^i}$ of $g$ to $D^i$ is a transversal metric on the foliation $\mathcal{F}_i$. The corresponding transversal Levi-Civita connection coincides with the connection on $D^i$ induced by $\nabla^0$, and its holonomy algebra coincides with $\h_i$.

Suppose that, for some $i$, it holds that $\h_i$ is a subalgebra of the unitary algebra $\u(D_x^i)$. Then the exists a $\nabla^{0}$-parallel field $J^i$ of complex structures on the distribution $D^i$, i.e., we obtain a K\"ahler transverse structure on the foliation $\mathcal{F}_i$. Let $K_i$ be the complex line bundle generated by the complex volume form $\Omega_i$ associated to $(g_{D^i},J^i)$. As we have seen in Section \ref{secTransStrFol}, the curvature of the induced connection on $K_i$ coincides with the Ricci form $\rho^i$, and the holonomy of this connection coincides with the projection of $\h_i$ to $\u(1)=\Real J^i_x$ with respect to the decomposition $\u(D_x^i)=\su(D_x^i)\oplus\u(1)$.


\begin{theorem}\label{thCriter} Let $(M,\theta,g)$ be a simply connected K-contact sub-Riemannian manifold with the complete Reeb vector field $\xi$. Then the following conditions are equivalent: 
\begin{itemize}
	\item[1.]	
	the horizontal holonomy algebra  $\hol(\nabla)$ is a codimension-one ideal of  $\hol(\nabla^0)$;
	
	\item[2.] it holds $\hol(\nabla^0)\subset\u(m)$, i.e., the Riemannian foliation $\mathcal{F}$ admits a transverse K\"ahler structure $J$; the space $D^0_x$ from decomposition \eqref{deRD} is trivial; each algebra $\h_i$ from  decomposition \eqref{deRh} is contained in $\u(D^i_x)$ and it contains  $\u(1)=\Real J_x^i$, where $J_x^i={J_x}|_{D^i_x}$;
	it holds  
	\begin{equation}\label{conddtheta=}d\theta=b_1\rho^1+\cdots+ b_r\rho^r,\quad b_1,\dots, b_r\in \Real\backslash\{0\}.\end{equation}
\end{itemize}	
\end{theorem}

\begin{rem}\label{rem1}
	If the conditions of Theorem \ref{thCriter} hold true, then the transverse K\"ahler structure $J$ is not defined uniquely. In fact, $\mathcal{F}$ admits $2^r$ K\"ahler structures, namely, each of them is of the form
	$$c_1J_x^1+\cdots+c_rJ_x^r,\quad c_1,\dots,c_r\in\{-1,1\}.$$ 
\end{rem}

\begin{theorem}\label{thholnabla} Let $(M,\theta,g)$ be a simply connected K-contact sub-Riemannian manifold with the complete Reeb vector field $\xi$. Suppose that the horizontal holonomy algebra  $\hol(\nabla)$ is a codimension-one ideal of  $\hol(\nabla^0)$. Let $J$ be a fixed transverse K\"ahler structure. Then there are decompositions \begin{align}\label{deRD1}  D_x&=D_x^1\oplus\cdots\oplus D_x^r,\\
\label{deRh1}\hol_x(\nabla^{0})&=\h_1\oplus\cdots\oplus \h_r,\\
	\label{deRh2}\hol_x(\nabla)&=\h'_1\oplus\cdots\oplus\h'_r\oplus\mathfrak{t}_x^\bot,\end{align} where, for each $i\in\{1,\dots,r\}$, 
	$$\h_i=\h'_i\oplus \Real J_x^i \subset\u(D^i_x),\quad J_x^i={J_x}|_{D^i_x},$$ is  an irreducible Riemannian holonomy algebra, and $\mathfrak{t}^\bot_x$ is the orthogonal complement to a one-dimensional subalgebra
	$$\t_x=\Real (a_1 J_x^1+\cdots +a_r J_x^r), \quad a_1,\dots,a_r\in\Real\backslash\{0\},$$
		in the commutative Lie algebra $\mathrm{span}\{J_x^1,\dots,J_x^r\}$.
	\end{theorem}	
 
\begin{rem}\label{rem2}
If $r\geq 2$, then the holonomy algebra $\hol(\nabla)$ from Theorem \ref{thholnabla} cannot be the holonomy algebra of a  Riemannian manifold. Indeed, $\hol(\nabla)$ preserves the decomposition \eqref{deRD1}, but, unlike $\hol(\nabla^0)$, $\hol(\nabla)$ is not the direct sum of subalgebras from $\so(D_x^i)$, i.e., $\hol(\nabla)$ does not satisfy the de Rham decomposition theorem. If  
$r=1$, then $\hol(\nabla^0)\subset\u(m)$ is an irreducible subalgebra. If $\hol(\nabla^0)\neq\u(m)$ is the holonomy algebra of a symmetric K\"ahler space, then $\hol(\nabla)=\su(m)\cap \hol(\nabla^0)$ is not a Riemannian holonomy algebra.
If $\hol(\nabla^0)=\u(m)$, then $\hol(\nabla)=\su(m)$, and this is the only case when the holonomy algebra $\hol(\nabla)$ from Theorem \ref{thholnabla} is a Riemannian holonomy algebra.	
\end{rem}

{\bf Proof of Theorems \ref{thCriter} and \ref{thholnabla}.}

Suppose that the horizontal holonomy algebra  $\hol_x(\nabla)$ is a codimension-one ideal of  $\hol_x(\nabla^0)$. Let $\t_x$ be as in Section \ref{secversus}.
It is obvious that $\t_x$ is contained in the direct sum of the centers of the Lie algebras $\h_i$.
If a subalgebra $\h_i\subset \so(D^i_x)$ has a non-trivial center, then $\h_i\subset \u(D_x^i)$, and the center of $\h_i$ coincides with $\u(1)=\Real J^i_x$. Thus, (after a change of the enumeration in \eqref{deRD}) there exists an $s$, $1\leq s\leq r$ such that
$$\t_x=\Real A_x\quad A_x= a_1 J_x^1+\cdots +a_s J_x^s, \quad a_1,\dots,a_s\in\Real\backslash\{0\}.$$
Let $$A=a_1 J^1+\cdots +a_s J^s.$$
Normalizing $A$, we may assume that $$g(A,A)=1.$$
It is clear that the line bundle
$$\t=\left<A\right>=\left<a_1 J^1+\cdots +a_s J^s\right>$$
is $\nabla^0$-parallel. 
Let $y\in M$ be an arbitrary point and let $\mu$ be a horizontal curve from $x$ to $y$. As it is well-known,
\begin{align*} \hol_y(\nabla)&=\tau_\mu(\hol_x(\nabla)),\\ 
\hol_y(\nabla^0)&=\tau_\mu^0(\hol_x(\nabla^0))=\tau_\mu(\hol_x(\nabla^0)).
\end{align*}
We conclude that 
\begin{equation}\label{hol_y}\hol_y(\nabla^0)=\hol_y(\nabla)\oplus\t_y.\end{equation}
From the results of Section \ref{secExtCon} it follows that
\begin{equation}\label{RWR0} R^W(\alpha)=R^0(\alpha)+d\theta(\alpha) N, \quad \forall \alpha\in\Gamma(\wedge^2 D),\end{equation}
where $N$ is the field of endomorphism defining the Wagner connection, and $R^W$ is the curvature tensor of the Wagner connection. 
From Theorem \ref{thholW}  it follows that
$$R_y^W(\alpha)\in\hol_y(\nabla),\quad\forall y\in M,\quad\forall \alpha\in\wedge^2D_y.$$
From this and \eqref{hol_y} it follows that
$$\pr_\t (R^0(\alpha))+d\theta(\alpha) \pr_\t N=0, \quad \forall \alpha\in\Gamma(\wedge^2 D),$$
where the projection is taken with respect to the decomposition \eqref{hol_y}.
Let $$K:\wedge^2\Real^{2k}\to\u(k)=\su(k)\oplus\Real J$$
be any algebraic curvature tensor, then it is easy to check that \begin{equation}\label{prJ} \pr_{\Real J}R(\alpha)=-\frac{1}{k}\rho(\alpha) J,\quad \forall\alpha\in\wedge^2\Real^{2k}.\end{equation}
This implies that
$$\pr_{\left<J^1\right>\oplus\cdots\oplus\left<J^s\right>} (R^0(\alpha))=-\frac{1}{m_1}\rho^1(\alpha)J^1-\cdots-\frac{1}{m_s}\rho^s(\alpha)J^s, \quad \alpha\in\Gamma(\wedge^2 D),$$
where the numbers $2m_i$ are the ranks of the distributions $D^i$. Using this, we get
\begin{multline*}\pr_\t (R^0(\alpha))=g\left(-\frac{1}{m_1}\rho^1(\alpha)J^1-\cdots-\frac{1}{m_s}\rho^s(\alpha)J^s ,A\right)A\\ 
=\left(-\frac{a_1}{m_1}\rho^1(\alpha)g(J^1,J^1)-\cdots-\frac{a_s}{m_s}\rho^s(\alpha)g(J^s,J^s)\right)A=-2
\left(a_1\rho^1(\alpha)+\cdots+a_s\rho^s(\alpha)\right)A.
\end{multline*}
On the other hand, $$\pr_\t N=2\varphi A, \quad\varphi =\frac{1}{2}g(N,A).$$
We conclude that 
$$\varphi d\theta=a_1\rho^1+\cdots+a_s\rho^s.$$
We claim that the function $\varphi$ is a non-zero constant. Indeed, applying the exterior derivative to the last equality, we get that $$d\varphi\wedge d\theta=0.$$
Since  $\iota_\xi d\theta=0$, we see that $\iota_\xi d\varphi=0$, i.e., $d\varphi\in\Gamma(D^*)$. Since the restriction of $d\theta$ to $\wedge^2 D$ is non-degenerate, 
and $2m=\mathrm{rk}D\geq 4$, we conclude that $d\varphi=0$, i.e., $\varphi$ is a constant.
 If $\varphi=0$, then $\pr_\t N=0$, $R^0$ take values in $\hol(\nabla)$, and we get a contradiction.  Since $d\theta\in\Gamma(\wedge^2 D)$ is non-degenerate, we conclude that $r=s$ and the distribution $D^0$ is trivial. Thus, the first condition of Theorem \ref{thCriter} implies the second condition.
 
 Suppose that the second condition of Theorem \ref{thCriter} holds true. Let 
 $$A=a_1 J^1+\cdots +a_r J^r$$ be the normalization of  
$$b_1 J^1+\cdots +b_r J^r$$
and consider the $\nabla^0$-parallel bundle $\t=\left<A\right>$. Let $\alpha\in\Gamma(\wedge^2(D))$.
Using \eqref{RWR0}, \eqref{NWagner} and the just used computations, we get
$$g(R^W(\alpha),A)=0.$$
From this and Theorem \ref{thASnabla} it follows that $\hol(\nabla)$ is orthogonal to $\t$. Since $d\theta$ is non-degenerate, each $\rho^i$ is non-zero. This  and \eqref{prJ} imply that $\hol(\nabla^0)$ contains $\Real J^i$. Hence 
$\hol(\nabla^0)$ contains $\t$, and $\hol(\nabla^0)\neq \hol(\nabla)$. From this and Theorem \ref{thHolnabl} it follows that $\hol(\nabla)\subset\hol(\nabla^0)$ is a codimension-one ideal. \qed

\begin{prop}\label{propholnabinsu} Under the conditions of Theorem \ref{thholnabla}, $\hol_x(\nabla)$ is contained in $\su(D_x)=\su(m)$ for a transverse K\"ahler structure $J$ if and only if $$d\theta=b\rho,\quad b\in\Real\backslash\{0\},$$ where $\rho$ is the Ricci form corresponding to $J$. 
	In this case $$\hol_x(\nabla^0)=\hol_x(\nabla)\oplus\Real J_x.$$
\end{prop}

{\bf Proof.} The condition $\hol_x(\nabla)\subset\su(m)$ is equivalent to the condition $\t_x^\bot\subset\su(m)$, which is equivalent to the condition $\t_x=\Real J_x$. \qed

\section{Consequences and examples}\label{secEx}
 In this section we give several examples and applications.

\subsection{A construction}

Let us give a construction of spaces satisfying the equivalent conditions of Theorem \ref{thCriter}. 
Let $(M_1,g_1)$,...,$(M_r,g_r)$ be K\"ahler manifolds such that their Ricci forms $\rho^i$ are non-degenerate and exact, i.e.,
$$\rho^i=d\theta^i,\quad \theta^i\in\Omega^1(M_i).$$
Consider the manifold
$$M=\Real\times M_1\times\cdots\times M_r.$$
Fix non-zero numbers $b_1,\dots,b_r$, define the 1-form
$$\theta=dt+b_1\theta^1+\cdots+b_r\theta^r\in\Omega^1(M),$$
where $t$ is the coordinate on $\Real$. It is obvious that $\theta$ is a contact form on $M$.
We consider the tangent bundles $TM_i$ as the subbundles of $TM$. Define the distributions 
$$D_i=\{X-b_i\theta(X)\partial_t|X\in TM_i\}\subset TM.$$
Each distribution $D_i$ may be identified with the tangent bundle $TM_i\subset TM$. Under this identification, each Riemannian metric $g_i$ defines the sub-Riemannian metric on the distribution $D^i$, and the corresponding horizontal connection on $D^i$ corresponds to the Levi-Civita connection on $(M_i,g_i)$. 
Finally, the contact distribution
$$D=D^1\oplus \cdots \oplus D^r=\ker\theta$$ with the sub-Riemannian metric $$g_1+\cdots+g_r$$ satisfies the second condition of Theorem \ref{thCriter}.

\subsection{Sasaki structures}

Recall that a Sasaki manifold $(M,\xi,h)$ is a Riemannian manifold $(M,h)$ with a unit Killing vector field $\xi$ such that the field of endomorphisms $\Phi$ defined by
$$\Phi(X)=\nabla^h_X\xi,\quad X\in\Gamma(TM),$$  satisfies the conditions
\begin{align*}
\Phi^2&=-\mathrm{id}_{TM}+\theta\otimes \xi,\\
(\nabla^h_X\Phi)Y&=h(\xi,Y)X-h(X,Y)\xi,\quad\forall X,Y\in\Gamma(TM),\\
\end{align*}    
where $\theta$ is the 1-form $h$-dual to the vector field $\xi$.

Let $(M,\theta,g)$ be a K-contact sub-Riemannian manifold. 
Consider the Riemannian metric
$$h=\theta\otimes\theta+g.$$ It is not hard to check the following relation between the Levi-Civita connection $\nabla^h$ and the Schouten connection $\nabla$:
\begin{align*}
\nabla^h_XY&=\nabla_XY-d \theta(X,Y)\xi,\\
\nabla^h_X\xi&=\psi(X),\\
\nabla^h_\xi X&=[\xi,X]+\psi(X),
\end{align*}
where $X,Y\in\Gamma(D)$, and the endomorphism $\psi$ is defined by the equality \begin{equation}\label{psi}g(\psi(X),Y)=d\theta(X,Y).\end{equation}
Using this, it is easy to prove the following:
\begin{prop} Let $(M,\theta,g)$ be a K-contact sub-Riemannian manifold with the Reeb vector field $\xi$. Then $(M,\xi,h)$ is a Sasaki manifold if and only if the endomorphism $\psi$ given by \eqref{psi} satisfies the conditions $$\psi^2=-\mathrm{id}_D,\quad \nabla\psi=0.$$
\end{prop}

Until the end of this section we suppose that $(M,\theta,g)$ has the structure of a Sasaki manifold with a complete Reeb vector field $\xi$. The holonomy algebra $\hol(\nabla^0)$ is contained in $\u(m)$, see \cite{HV}.

\begin{cor}\label{corS1} The condition $\hol(\nabla)\neq \hol(\nabla^0)$ holds true if and only if the transversal K\"ahler structure of each foliation $\mathcal{F}_i$ defined in Section \ref{secClassif} is Einstein and not Ricci-flat.
		\end{cor} 

{\bf Proof.} Suppose that $\hol(\nabla)\neq \hol(\nabla^0)$. From Remark \ref{rem1} it follows that the transverse K\"ahler structure $J$ from Theorem \ref{thCriter} may be taken to be the transverse K\"ahler structure defined by the Sasaki structure. Then it holds 
$$g(X,Y)=d\theta (JX,Y),\quad X,Y\in\Gamma(D).$$ 
From this and \eqref{conddtheta=} it follows that
$$g=b_1\mathrm{Ric}_1+\cdots+b_r\mathrm{Ric}_r,\quad b_1,\dots , b_r\in\Real\backslash\{0\}.$$
Conversely, the last equality implies \eqref{conddtheta=}, and by Theorem \ref{thCriter} it holds $\hol(\nabla)\neq \hol(\nabla^0)$. \qed

\begin{cor}\label{corS2}    
	The condition $\hol(\nabla)\subset\su(m)$ holds true if and only if the transversal K\"ahler structure of the characteristic foliation $\mathcal{F}$  is Einstein.
\end{cor} 

{\bf Proof} The condition $\hol(\nabla^0)\subset\su(m)$ holds true if and only if the transversal K\"ahler structure of  $\mathcal{F}$ is Ricci-flat \cite{HV}. This proves the corollary in the   case  $\hol(\nabla)= \hol(\nabla^0)$. If  
$\hol(\nabla)\neq \hol(\nabla^0)$, then the statement follows from Proposition \ref{propholnabinsu} and Corollary \ref{corS1}. \qed

 Let $A$ be one of the Lie groups $S^1$ or $\Real$. Recall that a Sasaki manifold $(M,\xi,h)$ is called regular if the Reeb vector field generates a free and proper action of a Lie group $A$. A regular Sasaki manifold $(M,\xi,h)$ is called {\it  a contactization} of a K\"ahler manifold $(Q,b,\omega)$, if there is a map  $p: M\to Q$ defining a  principle  $A$-bundle such that $$g=p^* b,\quad d\theta=p^*\omega,$$ see, e.g., \cite{ACHK}. The projection $p$ may be used to identify fibres of the contact distribution $D$  with the corresponding tangent spaces to $Q$. From \eqref{SchoutenCon} it follows that if $X,Y\in\Gamma(D)$, then $$p\nabla_XY=\nabla^b_{pX}(pY).$$ From \eqref{Schouten} it follows that
 $$p R(X,Y)Z=R^b(pX,pY)(pZ),\quad X,Y,Z\in\Gamma(D).$$
 This and the above corollaries imply
 
 \begin{cor}\label{corS3} Let a regular Sasaki manifold $(M,\xi,h)$ be  a contactization of a complete and simply connected K\"ahler manifold $(Q,b,\omega)$.
 	Then $\hol(\nabla)\neq \hol(\nabla^0)$ if and only if $(Q,b,\omega)$ is a product of K\"ahler-Einstein manifolds.
 	\end{cor}

\begin{cor}\label{corS4}  Let a regular Sasaki manifold $(M,\xi,h)$ be  a contactization of a complete and simply connected K\"ahler manifold $(Q,b,\omega)$.
	Then $\hol(\nabla)\subset\su(m)$ if and only if $(Q,b,\omega)$ is a  K\"ahler-Einstein manifold.
\end{cor}

\subsection{Parallel horizontal spinors}

Let $(M,\theta,g)$ be a simply connected contact sub-Riemannian manifold with the contact distribution $D$. The form $\theta$ defines an orientation on $D$. Let $P_{\SO}$ be the $\SO(2m)$-principle bundle of positively oriented frames on $D$.  A {\it spin-structure} on $D$ is a  reduction of $P_{\SO}$ to a $\Spin(2m)$-principle bundle $P_{\Spin}$ corresponding to the two-fold covering $$\lambda:\Spin(2m)\to\SO(2m).$$
Suppose that a spin-structure on $D$ is fixed.
Denote by $\Delta_{2m}$ the complex 
$\Spin(2m)$-spinor module. {\it The spinor bundle} $S$  is defined as the associated bundle,
$$S=P_{\Spin}\times_\lambda \Delta_{2m}.$$
The connections $\nabla$ and $\nabla^0$ are extended to the connections
$$\nabla^S:\Gamma(D)\times\Gamma(S)\to\Gamma(S),$$
$$\nabla^{0,S}:\Gamma(TM)\times\Gamma(S)\to\Gamma(S).$$
The holonomy algebras of these connections are respectively the representations of $\hol_x(\nabla)$ and 
$\hol_x(\nabla^0)$ on $S_x\cong\Delta_{2m}$. By the holonomy principle, there exists a non-zero spinor $s\in \Gamma(S)$ such that $\nabla^S s=0$ (resp., $\nabla^{0,S}s=0$) if and only if $\hol_x(\nabla)$ (resp., $\hol_x(\nabla^0)$) preserves a non-zero element $s_x\in S_x$.

From now on suppose that $(M,\theta,g)$ is a 
simply connected K-contact sub-Riemannian manifold.
As we have seen, $\hol_x(\nabla^0)\subset\so(2m)$ is a Riemannian holonomy algebra. Consequently, there exists a  non-zero $\nabla^{0,S}$-parallel spinor if and only if each Lie algebra in decomposition \eqref{deRh} is one of $\su(k)$, $\mathfrak{sp}(k)$, $\mathfrak{g}_2$, or $\mathfrak{spin}(7)$, see \cite{Wang}. If   $\hol_x(\nabla)=\hol_x(\nabla^0)$, then $\nabla^S$-parallel spinors are exactly  $\nabla^{0,S}$-parallel spinors. Consider the case when  
$\hol_x(\nabla)\neq\hol_x(\nabla^0)$.

\begin{theorem}\label{thparalspin} Let $(M,\theta,g)$ be a 
	simply connected K-contact sub-Riemannian manifold of dimension $2m+1\geq 5$ with a fixed spin structure on the contact distribution $D$ and the complete Reeb vector field $\xi$. Suppose that $\hol_x(\nabla)\neq\hol_x(\nabla^0)$. Then there exists a non-zero $\nabla$-parallel spinor if and only if there exists a transverse K\"ahler structure $J$ such that $\hol_x(\nabla)\subset\su(m)$.
\end{theorem}


{\bf Proof of Theorem \ref{thparalspin}.}

It is well-known that the $\u(n)$-module $\Delta_{2n}$ admits invariant decomposition 
\begin{equation}\label{Delta2n}\Delta_{2n}=\bigoplus_{k=0}^n\wedge^k\Co^n\end{equation}
such that $\Co^n$ is the standard $\su(n)$-module, and the complex structure $J\in\u(n)$ acts on $\wedge^k\Co^n$
as the multiplication by $(m-2k)i$.

\begin{lem} Let $\h=\h'\oplus\Real J\subset\u(n)$ be the holonomy algebra of an irreducible K\"ahler symmetric space, then $\h'$ does not annihilate any non-zero element from $\wedge^k\Co^n$ for $k=1,\dots,n-1$.
\end{lem}

{\bf Proof.} It is enough to show that $\h'$ does not annihilate any non-zero element from $\wedge^k(\Co^n)^*$, which is the same as to show that  $\ff=\h'\otimes \Co$
has this property. It is known that if $\ff\neq\mathfrak{sl}(n,\Co)$, then the first prolongation of the subalgebra $\ff\oplus\Co\id\subset\gl(n,\Co)$ contains a submodule $U$ isomorphic to $(\Co^n)^*$ such that $$\ff\oplus\Co\id=\{\varphi(X)|\varphi\in U,X\in\Co^n\},$$ see,  \cite[pages 95, 99]{M-Sch99}. 
Suppose that $\ff$ preserves a form $\omega\in \wedge^k(\Co^n)^*$, where $1\leq k\leq n-1$. Then the Lie algebra $\ff\oplus\Co\id$ preserves the line $\Co\omega\subset \wedge^k(\Co^n)^*.$ 
For each $\varphi\in U$, there exists $\lambda_\varphi\in (\Co^n)^*$ such that
$$(\varphi(X)\cdot\omega)(X_1,\dots,X_k)=
\lambda_\varphi(X)\omega(X_1,\dots,X_k),\quad\forall X,X_1,\dots,X_k\in\Co^n.$$
Alternating this equality and using
$$\varphi(X)Y=\varphi(Y)X,\quad\forall X,Y\in\Co^n,$$
we obtain 
$$\lambda_\varphi\wedge \omega=0.$$
It is clear that the map $$U\to(\Co^n)^*,\quad \varphi\mapsto\lambda_\varphi$$ is an isomorphism of  $\ff\oplus\Co\id$-modules. Thus, it holds that
$$\omega_1\wedge \omega=0,\quad\forall\omega_1\in (\Co^n)^*.$$
 Since $1\leq k\leq n-1$, this implies that $\omega=0$. \qed

Since $\hol_x(\nabla)\neq\hol_x(\nabla^0)$, the second condition of Theorem \ref{thCriter} holds true. Fix a transversal K\"ahler structure $J$. 
Decomposition \eqref{deRD1} implies the isomorphism
$$\Delta_{2m}=\Delta^1\otimes\cdots\otimes\Delta^r$$
of $\hol_x(\nabla)$-modules, 
where $\Delta^i$ is the spinor module corresponding to $D^i_x$, $i=1,\dots,r$.
Decomposition \eqref{Delta2n} then implies
$$\Delta_{2m}=\bigoplus_{k_1,\dots, k_r}\wedge^{k_1}(D_x^1\otimes \Co)\otimes\cdots\otimes\wedge^{k_r}(D_x^r\otimes \Co).$$
Fix a non-zero spinor $s$ from one of these submodules and suppose that $\hol_x(\nabla)$ annihilates $s$.
From decomposition \eqref{deRh2} and the lemma it follows that  $s$ belongs to a submodule corresponding to $k_i\in\{0,m_i\}$, where $2m_i=\dim D^i_x$, $i=1,\dots,r$.
As above, $\t_x=\Real A_x$, where 
$$A=a_1 J^1+\cdots +a_r J^r.$$
An element $$C=c_1{J^1_x}+\cdots+ c_r{J^r_x}$$
belongs to $\t^\bot$ if and only if $g(C,A_x)=0$, that is
\begin{equation}\label{eq1}
m_1a_1c_1+\cdots+m_ra_rc_r=0.
\end{equation}
The condition $Cs=0$ takes the form
\begin{equation}\label{eq2}
(m_1-2k_1)c_1+\cdots+(m_r-2k_r)c_r=0.
\end{equation}
Each solution of the equation \eqref{eq1} is a solution of 
the equation \eqref{eq2} if and only if
$$\frac{m_1-2k_1}{m_1a_1}=\cdots= \frac{m_r-2k_r}{m_ra_r}.$$
Since $k_i\in\{0,m_i\}$,  each number $m_i-2k_i$ is $\pm m_i$.
If the last equality holds, then we may change the K\"ahler structure $J$ in such a way that $$a_1=\dots =a_r=1.$$ This implies that $\t_x=\Real J_x$, and  $\hol_x(\nabla)\subset\su(m)$. Conversely, if  $\hol_x(\nabla)\subset\su(m)$, then $\hol_x(\nabla)$ annihilates a spinor $s\in\Delta_{2m}$ \cite{Wang}. \qed

\end{document}